\documentclass[12pt]{amsart}
\usepackage{amssymb}
\newtheorem{thm}{Theorem}[section]
\newtheorem{crl}[thm]{Corollary}
\newtheorem{prp}[thm]{Proposition}
\newtheorem{lm}[thm]{Lemma}

\newcommand{\der }{\partial }

\newcommand{\limline}{\limits_{-\!-\!-\!-\!-\!-}}
\newcommand{\limsim}{\limits_{\sim\!\sim\!\sim\!\sim\!\sim\!\sim}}
\newcommand{\limeq}{\limits_{=\!=\!=\!=\!=\!=}}
\newcommand{\limsimeq}{\limits_{\simeq\!\simeq\!\simeq\!\simeq\!\simeq\!\simeq}}
\newcommand{\limequiv}{\limits_{\equiv\!\equiv\!\equiv\!\equiv\!\equiv\!\equiv}}
\newcommand{\limsmile}{\limits_{\smile\!\smile\!\smile\!\smile\!\smile}}
\newcommand{\limapprox}{\limits_{\approx\!\approx\!\approx\!\approx\!\approx\!\approx}}
\newcommand{\limcong}{\limits_{\cong\!\cong\!\cong\!\cong\!\cong\!\cong}}
\newcommand{\limasymp}{\limits_{\asymp\!\asymp\!\asymp\!\asymp\!\asymp\!\asymp}}
\newcommand{\limdoteq}{\limits_{\doteq\!\doteq\!\doteq\!\doteq\!\doteq\!\doteq}}
\newcommand{\limfrown}{\limits_{\frown\!\frown\!\frown\!\frown\!\frown}}

\newcommand{\limli}{\limits_{-\;-\;-\;-}}

\begin{document}
\title{Wronskians as $n$-Lie multiplications}

\author{A.S. Dzhumadil'daev}
\address{Institute of Mathematics, Pushkin str.125, Almaty,
KAZAKHSTAN}
\email{askar@math.kz}
\maketitle

\begin{abstract}
Filipov proved that Jacobian algebra is $n$-Lie.
In our paper we consider algebras defined on associative commutative algebra $U$ with derivation $\der$ by $(k+1)$-multiplication
 $V^{0,1,\ldots,k}=\der^0\wedge\der^1\wedge\cdots\wedge \der^k$ (Wronskian).
We study whether they have $(k+1)$-Lie, $k$-left commutative and homotopical $(k+1)$- Lie structures. We prove that:
Wronskian algebra $(U,V^{0,1,\ldots,k})$
is homotopical  $(k+1)$-Lie, moreover,
$(U,\{0,\lambda_2V^{0,1},\ldots,\lambda_{i+1}V^{0,1,\ldots,i},\ldots ; i=1,2,\ldots\})$ is homotopical for any $\lambda_2,\lambda_3,\ldots \in K$;
$k$-left commutative, if
and only if $k>2$;  it is $(k+1)$-Lie, if and only if $k=1,$ $p$ is any or
 $(p,k)=(2,2)$ or $(p,k)=
(2,4)$ or $(p,k)=(3,3),$ where $p$ is the characteristic of main field $K$.
\end{abstract}

\section{Introduction}

Let $U$ be an associative commutative algebra with commuting derivations
$\der_1,\ldots,\der_n.$ The following determinant on $U,$ called Jacobian,
is well known
$$Jac_n (u_1,\ldots,u_n)=\left|
\begin{array}{ccc}
\der_1u_1&\cdots&\der_1u_n\\
\vdots&\vdots&\vdots\\
\der_nu_1&\cdots&\der_nu_n\\
\end{array}
\right|$$
In \cite{Filipov1}, \cite{Filipov2} it was proved that
$(U,Jac_n)$ as $n$-ary algebra is $n$-Lie, i.e.,
Jacobian satisfies the Leibniz rule:
$$Jac_n(u_1,\ldots,u_{n-1},Jac_n(u_n,\ldots,u_{2n-1}))=$$
$$\sum_{i=1}^nJac_n(u_n,\ldots,u_{n+i-1},Jac_n(u_1,\ldots,u_{n-1},u_{n+i}),u_{n+i+1},\ldots,u_{2n-1}).$$

There exists another remarkable determinant, called Wronskian.
It is defined on any associative commutative algebra $U$ with one
derivation $\der$ by the rule
$$V^{0,1,\ldots,k}(u_0,\ldots,u_k)=$$
$$\left|
\begin{array}{cccc}
u_0&u_1&\cdots&u_k\\
\der\,u_0&\der\,u_1&\cdots&\der\,u_k\\
\vdots&\vdots&\vdots&\vdots\\
\der^ku_0&\der^ku_1&\cdots&\der^ku_k\\
\end{array}\right|
$$
The aim of our paper is to study algebra $U$ as $n$-ary algebras
under wronskian multiplications.

Let $A$ and $M$ be vector spaces and $T^k(A,M)=Hom(A^{\otimes k},M)$ be a space of polylinear maps $A\times \cdots\times A\rightarrow A$ with $k$ arguments.
Set $T^0(A,M)=M$ and $T^k(A,M)=0,$ if $k<0.$ Set $T^*(A,M)=\oplus_kT^k(A,M).$

Let $\wedge^kA$ be $k$-th exterior power of $A$ and $C^k(A,M)=Hom(\wedge^kA,M)$ be subspace of $T^k(A,M).$ Set $C^0(A,M)=M$ and $C^k(A,M)=0,$ if $k<0$ and
$C^*(A,M)=\oplus_kC^k(A,M).$

Let $A$ be algebra with signature $\Omega$ \cite{Kurosh}. This means that
$\Omega$ is a set of polylinear maps $A\times \cdots\times A\rightarrow A.$
An element $\omega\in \Omega$ is called as {\it $n$-ary multiplication},
if $\omega\in T^n(A,M).$ We set $|\omega|=n,$ if $\omega\in T^n(A,A).$
If necessary to pay attention to signature, instead of $A$ we
will write $(A,\Omega).$ If $\Omega$ consists of one element $\omega$ we set
$A=(A,\omega).$

Let $(A,\omega)$ an $n$-ary algebra with vector space
$A$ over a field $K$ of characteristic $p\ge 0$ and
polylinear map with $n$ arguments $A\times \cdots \times A\rightarrow A.$ Recall that a linear map $D: A\rightarrow A$ is called derivation of $A,$ if $$D(\omega(a_1,\ldots,a_n))=\sum_{i=1}^n\omega(a_1,\ldots,a_{i-1},D(a_i),a_{i+1},\ldots, a_n),$$
for any $a_1,\ldots,a_n\in A.$ Let $L_{a_1\ldots a_{n-1}}: A\rightarrow A$ be a linear map defined by $$L_{a_1,\ldots,a_{n-1}}a=\omega(a_1,\ldots,a_{n-1},a).$$
Let $Der\, A$ be an algebra of all derivations of the algebra $(A,\omega).$
An algebra $(A,\omega)$ is called {\it $n$-Lie} \cite{Filipov1}, if
$\omega$ is skew-symmetric and $$L_{a_1,\ldots,a_{n-1}}\in Der\,A,$$
for any $a_1,\ldots,a_{n-1}\in A.$ About $n$-Lie algebras called also as
Nambu, Filipov, Takhtajan see also \cite{Nambu}, \cite{Takh}.

Let $Sym_k$ be the permutation group and $sign\,\sigma$ is signum of the permutation $\sigma\in Sym_k.$ Denote by $Sym_{k,l}$ the subset of $k,l$-shufle
permutations, i.e., permutations $\sigma\in Sym_{k+l},$ such that $\sigma(1)<\cdots<\sigma(k), \sigma(k+1)<\cdots<\sigma(k+l).$
Usually by the set of permutations for $Sym_k$ one
understands the standard set $\{1,\ldots,k\},$ but in our paper we use some
non-standard sets of order $k,$ something like $2,3,\ldots,k+1\}$ in definition of $Q_l.$
Form the context will be clear what kind of sets we use.

An algebra $(A,\omega)$ with $n$-ary multiplication $\omega$ is called {\it $(n-1)$-left commutative}, if
$$\sum_{\sigma\in Sym_{2n-2}}sign\,\sigma\,\omega(a_{\sigma(1)},\ldots,a_{\sigma(n-1)},\omega(a_{\sigma(n)},\ldots,a_{\sigma(2n-2)},a_{2n-1}))=0,$$
for any $a_1,\ldots,a_{2n-2}, a_{2n-1}\in A.$

Call an algebra $(A,\omega)$ as {\it homotopical $n$-Lie} \cite{StashefLada}, if $\omega\in C^n(A,A)$ and
$$\sum_{\sigma\in Sym_{n-1,n-1}}sign\,\sigma\,\omega(a_{\sigma(1)},
\ldots,a_{\sigma(n-1)},\omega(a_{\sigma(n)},\ldots,a_{\sigma(2n-2)},
a_{\sigma(2n-1)}))=0,$$
for any $a_1,\ldots,a_{2n-2}, a_{2n-1}\in A.$
In \cite{Hanlon} such algebra is called Lie $n$-algebra.

An algebra $(A,\Omega)$ is called {\it homotopical} Lie algebra
\cite{StashefLada},
if $\Omega$ consists of elements $\omega_1,\omega_2,\ldots,$ such that
$\omega_k\in C^k(A,A)$  and
$$\sum_{\sigma\in Sym_{i-1,j}}sign\,\sigma\,\omega_i(a_{\sigma(1)},\ldots,a_{\sigma(i-1)},
\omega_j(a_{\sigma(i)},\ldots,a_{\sigma(i+j-1)}))=0,$$
for any $a_1,\ldots,a_{i+j-1}\in A$ and $i,j=1,2,\ldots .$

In section \ref{nlienleft} we prove that any $n$-Lie algebra is
$(n-1)$-left commutative and any $(n-1)$-left commutative algebra is
homotopical $n$-Lie. Examples of Wronskian algebras show that inverse of
these statements are not true.

Let $U$ be associative commutative algebra with derivation $\der.$
Let $V^{i_1,\ldots,i_k}=\der^{i_1}\wedge\cdots\wedge\der^{i_k}$ be
the {\it generalized Wronskian,} i.e.,
$$V^{i_1,\ldots,i_k}(u_1,\ldots,u_k)=\left|
\begin{array}{ccc}
\der^{i_1}u_1&\cdots&\der^{i_k}u_k\\ \vdots&\vdots&\vdots\\
\der^{i_k}u_1&\cdots&\der^{i_k}u_k\\
\end{array}\right|
$$
For example, $V^{0,1,2,\ldots,k}$ is the standard Wronskian.

\begin{thm} \label{trieste1} Let $U$ be any associative commutative algebra over the field $K$ of characteristic $p\ge 0$ and derivation
$\der.$ Then
\begin{enumerate}
\renewcommand{\theenumi}{\roman{enumi}}
\item For any $k>0,$ the algebra $(U,V^{0,1,\ldots,k})$ is homotopical $(k+1)$-Lie. Moreover, $(U,\{0,\lambda_{i+1}V^{0,1,\ldots,i}, i=1,2,\ldots\})$ is homotopical Lie, for any $\lambda_i\in k, \lambda_1=0.$
\item The algebra $(U,V^{0,1,\ldots,k}), k>0,$ is $k$-left
commutative, if and only if $k\ne 2.$
\item The algebra $(U,V^{0,1,\ldots,k})$ is $(k+1)$-Lie, if and only if
one of the following conditions are hold
\begin{itemize}
\item $k=1$ and  $p$ is any prime number or $0$
\item $k=2$ and $p=2$
\item $k=3$ and $p=3$
\item $k=4$ and $p=2.$
\end{itemize}
\end{enumerate}
\end{thm}

{\bf Remark.} One can check that $(U,\{V^{0,1}, \lambda_3V^{0,2,3}, \lambda_4V^{0,2,3,4}\})$ is homotopical, if and only if $-5 \lambda_3^2+ 7 \lambda_4=0$ and there are no homotopical prolongation for $\{V^{0,1}, V^{0,1,4}\}.$
Description  of homotopical structures with beginning part $V^{0,1}$ will be discussed in another paper. Results will follow from the following general fact.
Suppose that $A$ is $\Omega$-homotopical Lie algebra, where
$\Omega=\{\omega_k\in \wedge^k(A,A): k=1,2,\ldots\},$ such that
$\omega_1=0.$ Denote by $A^{lie}$ an algebra $A$ under Lie
multiplication $\omega_2.$ Then $\omega_k\in
Z^k(A^{lie},A^{lie}),$ for any $k\ge 2.$

By theorem \ref{trieste1}, Wronskians as $n$-Lie multiplications
in sense of Filipov, for $n>2$  can appear only in cases of small
characteristics $p=2,3.$ In \cite{Filipov1} it was established the
following result: if $A$ is $n$-Lie with multiplication $\omega$
then $A$ is $(n-1)$-Lie under multiplication $i(a)\omega,$ for any
$a\in A.$ Using this construction from $(k+1)$-Lie algebras
$V^{0,1,\ldots,k}$ one can obtain another $n$-Lie algebras with
$n\le k.$

\begin{thm}\label{trieste2} Let $n\ge 2$ and $p=char\,K\ge 0.$
The following generalized Wronskians are $n$-Lie multiplications
\begin{itemize}
\item $n=2$
\begin{itemize}
\item $p=2, V^{2^r-2^l,2^r}, 0\le l\le r$
\item $p=2, \sum_{i=1}^{2^l}V^{i,2^l+1-i}, 0<l$
\item $p=3, V^{2\cdot 3^r,3^{r+1}}, 0\le r$
\end{itemize}
\item $n=3$
\begin{itemize}
\item $p=2, V^{1,2,4}$
\item $p=2, V^{2,3,4}$
\item $p=2, \sum_{i=1}^{2^l}V^{0,i,2^l+1-i}, 0<l$
\item $p=3, V^{1,2,3}$
\end{itemize}
\item $n=4$
\begin{itemize}
\item $p=2, V^{1,2,3,4}$
\item $p=3, V^{0,1,2,3}$
\end{itemize}
\item $n=5, p=2, V^{0,1,2,3,4}$
\end{itemize}
Over a fields of characteristic $0$ Wronskian $V^{0,1,\ldots,k}$ can
serve as $n$-Lie multiplication only in case of Lie algebras, $n=2, k=1.$
\end{thm}

\begin{crl} Let $U$ be associative commutative algebra with derivation $\der$ 
over a field of characteristic $p=3.$ Endow $U$ by 
$3$-multiplication $V^{1,2,3}.$ Then $U$ is a triple Lie system.  
\end{crl}

It is well known that, if $D_1, D_2$ are derivations of some associative 
commutative algebra $U,$ then $D_1\wedge D_2: (u,v)\mapsto D_1(a)D_2(b)
-D_2(a)D_1(b)$ is a Lie multiplication. In case of small characteristics 
might happen that $D\wedge F$ is also Lie multiplication, if $D$ or $F$ 
is not a derivation.
 
\begin{crl} Let  $D\in Der\,U$ be a derivation of some associative 
commutative algebra $U.$  Then the multiplication $D^2\wedge D^3: (u,v)\mapsto 
D^2(u)D^3(v)-D^2(v)D^3(u)$ endows $U$ by a structure of Lie algebra, if $p=3.$
The multiplication $D^{2^k}\wedge D^{2^k-1}: (u,v)\mapsto 
D^{2^k}(u)D^{2^k-1}(v)-D^{2^k-1}(u)D^{2^k}(v)$ is also endows $U$ 
by a structure of Lie algebras, if $p=2.$
\end{crl}

Our calculations are based on two ideas. The first idea (a polynomial trick)
means the following. Suppose that some statement $\mathcal X$  about associative commutative algebra with unit $U$ and commuting derivations ${\mathcal D}=<\der_1,\ldots,\der_n>$ was obtained by using
\begin{itemize}
\item linear properties of $U$
\item associativity, commutativity and unit properties of $U$
\item Leibniz rule for derivations $\der_1,\ldots,\der_n.$
\item commuting property of derivations $[\der_i,\der_j]=0, i,j=1,2,\ldots,n$
\end{itemize}
Then this statement is true for any associative commutative algebra with unit and commuting derivations.

In particular, we can take as $U$ the polynomial algebra $K[x_1,\ldots,x_n]$
and $\der_i=\der/\der\,x_i$ or algebra of divided power polynomials
(in case of $p>0$)
$$O_n({\bf m})=\{x^{\alpha}=\prod_{i=1}^nx^{(\alpha_i)}: 0\le \alpha_i<p^{m_i}, {\bf m}=(m_1,\ldots,m_n)\},$$
$$x^{\alpha}x^{\beta}=\prod_{i=1}^n{\alpha_i+\beta_i\choose \alpha_i}x^{\alpha+\beta},$$
with special derivations
$$\der_i: x^{\alpha}\mapsto x^{\alpha-\epsilon_i}, \epsilon_i=(0,\cdots,0,\mathop{1}\limits_i,0,\cdots,0), i=1,\ldots,n.$$
We can consider the statement $\mathcal X$ for (divided) polynomial algebras.

The second idea concerns ${\mathcal D}$-invariant polynomials
\cite{Dzhumavestnik}. Let $\psi\in T^k(U,U)$  be polylinar map. It is called ${\mathcal D}$-invariant, if
$$\der\psi(u_1,\ldots,u_k)=\sum_{i=1}^k\psi(u_1,\ldots,u_{i-1},\der u_i,u_{i+1},\ldots,u_k),$$
for any $u_1,\ldots,u_k\in U$ and $\der\in {\mathcal D}.$ In other words, $\psi$ is ${\mathcal D}$-invariant, if any $\der\in {\mathcal D}$ is derivation for $\psi.$ Notice that $U$ is $\mathcal D$-graded:
$$U=\oplus_{s\ge 0} U_s, U_sU_l\subseteq U_{s+l},$$
$$U_0=<1>,$$
$$\der_iU_s\subseteq U_{s-1},$$
$$U^{\mathcal D}=\{u\in U: \der_iu=0, \forall i=1,\ldots,n\}=U_0.$$
For graded ${\mathcal D}$-invariant polylinear map $\psi\in T^k(U,U)$
denote by $\pi\psi$ a polylinear form $\pi\psi\in T^k(U,U_0),$ defined on
homogeneous basic elements $e_1,\ldots,e_k\in U$ by
$$\pi\psi(e_1,\ldots,e_k)=\psi(e_1,\ldots,e_k),$$
if $\psi(e_1,\ldots,e_k)\in U_0$ and
$$\pi\psi(e_1,\ldots,e_k)=0,$$
if $\psi(e_1,\ldots,e_k)\in \oplus_{s>0}U_s.$
Call $\pi\psi$ as a {\it support} of $\psi$ and $k$-typle of basic
homogeneous elements $(e_1,\ldots,e_k),$ such that
$\pi\psi(e_1,\ldots,e_k)\ne 0$ as a {\it supporting chain}.
Let $\Gamma$  be a set of supporting chains. Then \cite{Dzhumavestnik}
$\psi$ can be restored by support $\pi\psi$ in a unique way. Namely,
$$\psi(u_1,\ldots,u_k)=\sum_{\{\alpha_1,\alpha_2,\ldots,\alpha_k\}\in \Gamma}
\frac{\der^{\alpha_1}(u_1)}{\alpha_1!}\ldots\frac{\der^{\alpha_k}(u_k)}{\alpha_k!}\pi\psi(x^{\alpha_1},\ldots,x^{\alpha_k}).$$
So, to find $\mathcal D$-polylinear form it is enough to calculate
its support. We us this method in calculations of $Q\psi,$ $Q_{short}\psi,$
$Q_{long}\psi$ and $Q_{alt}\psi.$ In our case $n=1$ and the restoring formula looks more simpler
$$\psi=\sum_{i_1,\ldots,i_k\in \Gamma}\lambda_{i_1,\ldots,i_k}
\frac{\der^{i_1}u_1}{i_1!}\ldots\frac{\der^{i_k}u_k}{i_k!},$$
where $\lambda_{i_1,\ldots,i_k}=\pi\psi(x^{i_1},\ldots,x^{i_k})\in K.$
In case of divided polynomials here $x^{i}$ should be changed to $x^{(i)}$
and $\frac{\der^{i}u}{i!}$ to $\der^iu.$
Calculations of support chains and support forms was made by Matematica.

\section{\label{nlienleft} Connections between $n$-Lie, $(n-1)$-left commutative and homotopical $n$-Lie structures}

Define quadratic maps
$$Q, Q_{short}, Q_{long}, Q_{alt}: C^k(A,A)\rightarrow T^{2k-1}(A,A),$$
by
$$Q\psi(a_1,\ldots,a_{2k-1})=
\psi(a_1,\ldots,a_{k-1},\psi(a_k,\ldots,a_{2k-1}))$$
$$-\sum_{i=1}^k
\psi(a_{k},\ldots,a_{k+i-1},\psi(a_1,\ldots,a_{k-1},a_{k+i}),
a_{k+i+1},\ldots,a_{2k-1}),$$

$$Q_{long}\psi(a_1,\ldots,a_{2k-1})=$$
$$\sum_{\sigma\in Sym_{k-1,k}, \sigma(k-1)=2k-1}
sign\,\sigma\,\psi(a_{\sigma(1)},\ldots,a_{\sigma(k-2}),a_{\sigma(k-1)},
\psi(a_{\sigma(k)},\ldots,a_{\sigma(2k-1)})),$$

$$Q_{short}\psi(a_1,\ldots,a_{2k-1})=$$
$$\sum_{\sigma\in Sym_{k-1,k},\sigma(2k-1)=2k-1}
sign\,\sigma\,\psi(a_{\sigma(1)},\ldots,a_{\sigma(k-1}),\psi(a_{\sigma(k)},\ldots,a_{\sigma(2k-2)}, a_{2k-1})),$$

$$Q_{alt}\psi(a_1,\ldots,a_{2k-1})=$$
$$\sum_{\sigma\in Sym_{k-1,k}}
sign\,\sigma\,\psi(a_{\sigma(1)},\ldots,a_{\sigma(k-1)},\psi(a_{\sigma(k)},\ldots,a_{\sigma(2k-2)}, a_{\sigma(2k-1)})),$$

These definitions divide elements $\{a_1,\ldots,a_{2k-1}\}$ into two types.
Call elements of the $(k-1)$ elements subset $\{a_1,\ldots,a_{k-1}\}$ as a
{\it short} and elements
of the  $k$ elements subset
$\{a_k,\ldots,a_{2k-1}\}$ as a {\it long.}

Notice that
$Q\psi\in T^{2k-1}(A,A)$ is skew-symmetric by
short arguments and by long arguments,
i.e., by first $(k-1)$ arguments and last $k$ arguments,
if $\psi\in C^k(A,A).$
It is easy to see that, $Q_{short}\psi\in T^{2k-1}(A,A)$ is skew-symmetric by all short  arguments and skew-symmetric by all long arguments, except one, if $\psi\in C^k(A,A).$
Similarly, $Q_{long}\psi\in T^{2k-1}(A,A)$ is skew-symmetric by all long arguments and by all short arguments except one, if $\psi\in C^k(A,A).$
Notice that $Q_{alt}\psi\in C^{2k-1}(A,A).$

\begin{prp} \label{trieste3}
Suppose that $\omega\in C^k(A,A)$ and
one of the following conditions are held
\begin{itemize}
\item $p=0,k>2,$
\item $p>0, k\not\equiv 0,1 (mod\,p), k\equiv
1(mod\,2),$
\item $p>0, k\not\equiv -1,2(mod\,p), k\equiv
0(mod\,2).$
\end{itemize}
Take place the following statements.

i) If $Q\omega=0,$ then $Q_{short}\omega=0,$ and $Q_{long}\omega=0.$

ii) If $Q_{short}\omega=0$ or $Q_{short}\omega=0,$ then $Q_{alt}\omega=0.$
\end{prp}

{\bf Proof.} Call $\sigma\in Sym_{k-1,k}$ as a short permutation,
if $\sigma(k-1)=2k-1.$ To any short permutation $\sigma$
corresponds $k-2$-typle $r'(\sigma)$ and $k$-typle $r''(\sigma)$
defined  by $$r'(\sigma)=\{\sigma(1),\ldots,\sigma(k-2)\},$$
$$r''(\sigma)=\{\sigma(k),\ldots,\sigma(2k-1)\}.$$

Call $\sigma\in Sym_{k-1,k}$ as a long permutation, if
$\sigma(2k-1)=2k-1.$ To any long permutation $\sigma$ corresponds
$k-1$-typle $r'(\sigma)$ and $k-1$-typle $r''(\sigma)$ defined by
$$r'(\sigma)=\{\sigma(1),\ldots,\sigma(k-1)\},$$
$$r''(\sigma)=\{\sigma(k),\ldots,\sigma(2k-2)\}.$$

Notice that $$r'(\sigma)\cup r''(\sigma)=\{1,\ldots,2k-2\},$$ for
any $\sigma\in Sym_{k-1,k}.$ In other words, any $r'(\sigma)$ is
uniquely defined by $r''(\sigma).$

Call an element of the form
$$A_{\sigma}:=\psi(a_{\sigma(1)},\ldots,a_{\sigma({k-1})},
\psi(a_{\sigma({k})},\ldots, a_{\sigma(2k-1)}))$$ as a short
(long) $\sigma$-element, or just short (long) element,
if $\sigma$ is short (long) permutation. Let
$Sym_{k-1,k}^s$ be the set of short permutations and
$Sym_{k-1,k}^l$ be the set of long permutations. It is evident
that $$Sym_{k-1,k}=Sym_{k-1,k}^s\cup Sym_{k-1,k}^l$$ and
\begin{equation}\label{a}
Q\psi=Q_{short}\psi+Q_{long}\psi,
\end{equation}

The identity $Q\psi=0$ for $\psi\in C^k(A,A)$ gives us that
$$\psi(a_{\sigma(1)},\ldots,a_{\sigma({k-1})},\psi(a_{\sigma({k})},\ldots,
a_{\sigma({2k-1})}))=$$
\begin{equation}\label{d}\sum_{i=0}^{k-1}(-1)^{k-i-1}
\psi(a_{\sigma(k)},\ldots,\widehat{a_{\sigma(k+i)}},\ldots,a_{\sigma(2k-1)},
\psi(a_{\sigma(1)},\ldots,a_{\sigma(k-1)},a_{\sigma(k+i)})).\end{equation}

In particular, (\ref{d}) means that, any short $\sigma$-element can be
presented as a sum of $k$ long elements. More exactly, the short
element $A_{\sigma}$  is a sum of $k$ long elements of the form
$A_{\tau},$ where the long permutation $\tau$ satisfies the
condition $r'(\tau)\subset r''(\sigma).$ Since
$|r''(\sigma)\setminus r'(\tau)|=1,$ there exists exactly one
element, say $i,$ such that $r''(\sigma)=r'(\tau)\cup \{i\}.$ Then
$i\le 2k-2$ and $i$ is not equal to $k-1$ elements of $r'(\tau).$
So, there are $k-1$ possibilities to choose $i.$

In other words, according (\ref{d}) the element
$$Q_{long}\psi(a_1,\ldots,a_{2k-1})=$$
$$
\sum_{\sigma\in Sym_{k-1,k}^s}sign\,\sigma\,
\psi(a_{\sigma(1)},\ldots,a_{\sigma({k-1})},\psi(a_{\sigma({k})},\ldots,
a_{\sigma({2k-1})}))$$ can be presented in the form
$$(k-1)\sum_{\sigma\in Sym_{k-1,k}^l}
\psi(a_{\sigma(1)},\ldots,a_{\sigma({k-1})},\psi(a_{\sigma({k})},\ldots,
a_{\sigma({2k-1})}))$$ So, from conditions $\psi\in C^k(A,A),
Q\psi=0$ one can obtain that
\begin{equation}\label{may11}
Q_{long}\psi= (k-1)Q_{short}\psi.
\end{equation}

Let us give another interpretation of (\ref{d}). If $\sigma$ is a
long permutation, then $A_{\sigma}$ is a sum of $k-1$ short
elements and one long element $A_{\tilde\sigma},$ where
$$\tilde\sigma= \left(\begin{array}{cccccccc}
1&\cdots&k-1&k&\cdots&2k-2&2k-1\\
\sigma(k)&\cdots&\sigma(2k-2)&\sigma(1)&\cdots&\sigma(k-1)&2k-1\\
\end{array}\right)$$
So, $A_{\sigma}$ can be presented as a sum of short elements of
the form $A_{\tau},$ where $r'(\tau)\subset r''(\sigma).$ More
exactly, $r''(\sigma)\setminus r'(\tau)=\{i\},$ for some $i\in
\{1,2,\ldots,2k-2\}$ and $i$ can not be equal to one of $k-2$
elements of $r'(\tau).$ So, there are $k$ possibilities for the
long permutation $\sigma,$ such that $A_{\tau}$ can be one of
summands    of  $A_{\sigma}.$ Notice that
$$sign\,\sigma=(-1)^{k-1}sign\,\tilde\sigma.$$

Therefore, summation
of $sign\,\sigma\,A_{\sigma}$ for all $\sigma\in Sym_{k-1,k}^l$
according (\ref{d})  gives us that
\begin{equation}\label{may12}
Q_{short}\psi=k\,Q_{long}\psi+(-1)^{k-1}Q_{short}\psi.
\end{equation}
We see that the determinant of the system of linear equations
(\ref{may11}) and (\ref{may12}) is equal to
$$\left|\begin{array}{cc} 1&-k+1\\ k&-1-(-1)^k\\
\end{array}\right|=-k^2+k+1+(-1)^k$$
Therefore, from the conditions $Q\psi=0$ and $\psi\in C^k(A,A)$
imply the following identities $$Q_{long}\psi-2 Q_{short}\psi=0,\quad k\equiv
-1 (mod\, p), k\equiv 0(mod 2), p>0,$$ $$Q_{long}\psi+Q_{short}\psi=0, \quad
k\equiv 2(mod\,p), k\equiv 0(mod \,2), p>0,$$ $$Q_{long}\psi-
Q_{short}\psi=0,\quad k\equiv 0(mod\, p), k\equiv 1(mod 2), p>0,$$
$$Q_{long}\psi=0, \quad k\equiv 1(mod\,p), k\equiv 1(mod \,2), p>0,$$
and $$Q_{long}\psi=0, Q_{short}\psi=0,\quad p=0, k>2.$$ So, according
(\ref{a}) $$Q_{alt}\psi=0,$$ if $p=0,k>2$ or  $k\not\equiv
0,1(mod\,p), k\equiv 1(mod\,2),$ or $k\not\equiv -1,2(mod\,p),
k\equiv 0(mod\,2).$

\begin{crl} If an algebra $(A,\omega)$ is $n$-Lie, then it is
$(n-1)$-left commutative. If algebra $(A,\omega)$ $(n-1)$-left commutative,
then it is homotopical $n$-Lie.
\end{crl}
In particular, any $n$-Lie algebra is
homotopical $n$-Lie. As show theorems \ref{trieste1} and
\ref{trieste2} the inverse of this statement is not true.
Proposition \ref{trieste3} for 3-algebras was also noticed in \cite{Awata}, \cite{Hoppe}.

\section{Proof of theorem \ref{trieste1} and \ref{trieste2}}

Let $${\bf Z}^n=\{\alpha=(\alpha_1,\ldots,\alpha_n): \alpha_i\in {\bf Z}\}$$
and
$${\bf Z}^n_+=\{\alpha\in {\bf Z}^n: \alpha_i\ge 0, i=1,\ldots,n\}.$$
Set for $\alpha\in {\bf Z}^n,$
$$|\alpha|=\sum_{i=1}^n\alpha_i.$$

Let $U$ be associative commutative algebra with derivation $\der.$
Let $\psi \in C^k(U,U)$ and $g\in C^l(U,U).$ Define
$f\smile g, f\wedge g\in C^{k+l}(U,U)$ and bilinear maps $Q(f,g),Q_{alt}(f, g)\in C^{k+l-1}(U,U),
Q_{short}(f,g)\in T^{k+l-1}(U,U)$ by
$$f\smile g(u_1,\ldots, u_{k+l})=f(u_1,\ldots,u_k)g(u_{k+1},\ldots,u_{k+l}),$$

$$f\wedge g(u_1,\ldots,u_{k+l})=\sum_{\sigma\in Sym_{k+l}}sign\,\sigma\,
(f\smile g)(u_{\sigma(1)},\ldots,u_{\sigma(k+l)}).$$

$$f\star g(u_1,\ldots,u_{k+l-1})=
f(u_1,\ldots,u_{k-1},g(u_k,\ldots,u_{k+l-1})),$$

$$Q(f,g)=f(u_1,\ldots,u_{k-1},g(u_k,\ldots,u_{k+l-1}))-$$
$$\sum_{i=1}^lg(u_k,\ldots,u_{k+i-2},f(u_1,\ldots,u_{k-1},u_{k+i-1}),u_{k+i},\ldots,u_{k+l-1}),$$

$$Q_{alt}(f, g)(u_1,\ldots,u_{k+l-1})=
\sum_{\sigma\in Sym_{k-1,l}}
sign\,\sigma\,(f\star g)(u_{\sigma(1)},\ldots,u_{\sigma(k+l-1)}),$$

$$Q_{short}(f, g)(u_1,\ldots,u_{k+l-1})=
\sum_{\sigma\in Sym_{k-1,l-1}}
sign\,\sigma\,(f\star g)(u_{\sigma(1)},\ldots,u_{\sigma(k+l-2)}, u_{k+l-1}).$$

Notice that  definitions of $Q$ as  bilinear maps are compatible with definitions of $Q$ as quadratic maps in the previous section:
$Q(f,f)=Q(f), Q_{alt}(f,f)=Q_{alt}(f), Q_{short}(f,f)=Q_{short}(f).$

Since $(C^*(U,U),\cup)$ is associative, $(C^*(U,U),\wedge)$ is also
associative. Let
$$C^k_{loc,s}(U)=\{\der_{i_1}\wedge \cdots \wedge \der^{i_k}:
0\le i_1<\ldots <i_k, i_1+\cdots+i_k=s\}$$
and $$C^k_{loc}(U)=\oplus_sC^k_{loc,s}(U).$$
Notice that $V^{\alpha}\in C^k_{loc,|\alpha|}(U),$ for any $\alpha\in {\bf Z}^n_+.$
Set $|\psi|=s,$ if $\psi\in C^k_{loc,s}(U).$

Let $\psi\in T^k(A,A).$ Define $i(a)\psi\in T^{k-1}(A,A)$ by
$$i(a)\psi(a_1,\ldots,a_{k-1})=\psi(a,a_1,\ldots,a_{k-1}).$$

\begin{prp}\label{filipov1} If $\psi$ is $k$-Lie,  then
$i(a)\psi$ is $(k-1)$-Lie for any $a\in A.$
\end{prp}

{\bf Proof.} See \cite{Filipov1}.

Proposition \ref{filipov1} can be modified by the following way
\begin{lm}\label{17july} If $\psi$ is $k$-Lie, then for $(k-l)$-Lie multiplications
$\psi_l:=i(a_1)i(a_2)\cdots i(a_l)\psi$ take place the following
relations
$$Q(\psi_i,\psi_j)=0, \quad i\le j.$$
\end{lm}

\begin{lm}\label{altynbalatugankun} $C^k_{loc,s}(U)=0,$ if $s<k(k-1)/2.$
\end{lm}
{\bf Proof.}
If $0\ne \der^{i_1}\wedge\cdots\wedge\der^{i_k}\in C^k_{loc,s}(U),$
then $s=i_1+\ldots+i_k\ge 0+1+2+\cdots+ (k-1)=(k-1)k/2.$

\begin{lm} \cite{dzhumafunctanalis}  For any $\alpha\in {\bf Z}^k_+,$ and
$\beta\in {\bf Z}^l_+,$
$$Q_{alt}(V^{\alpha},V^{\beta})\in C^{k+l-1}_{loc,|\alpha|+|\beta|}(U).$$
\end{lm}

\begin{crl}\label{altynaitugankun} $$Q_{alt}(V^{0,1,\ldots,k},V^{0,1\ldots,l})=0,$$
for any $k,l>0.$
\end{crl}

{\bf Proof.} Notice that $|V^{0,1,\ldots,k}|=k(k+1)/2.$ Therefore
$$Q_{alt}(V^{0,1,\ldots,k},V^{0,1,\ldots,l})\in C^{k+l+1}_{loc,(k^2+l^2+k+l)/2}(U).$$
It is evident, that $(k^2+l^2+k+l)/2<(k+l+1)(k+l)/2.$
Therefore, by lemma \ref{altynbalatugankun},
$C^{k+l+1}_{loc, (k^2+l^2+k+l)/2}(U)=0.$ Thus,
$Q_{alt}(V^{0,1,\ldots,k},V^{0,1,\ldots,l})=0.$

\begin{lm}
\label{altynkustugankun} $Q_{short}(V^{0,1,\ldots,k},V^{0,1,\ldots,k})=0,$ if $k> 2.$ If $k=2,$ then $Q_{short}(V^{0,1,2},V^{0,1,2})=V^{0,1,2,3}\smile id.$
\end{lm}

{\bf Proof.} Notice that $Q_{short}(V^{0,1,\ldots,k},V^{0,1,\ldots,k})$ is
a linear combination of cochains of the form
$(\der^{i_1}\wedge \der^{i_{2k}})\smile \der^{i_{2k+1}},$
such that $i_1+\cdots+i_{2k}+i_{2k+1}=
k^2+k$ and $0\le i_1<i_2<\cdots <i_{2k}.$
We have $i_1+\cdots+i_{2k}>0+1+2+\cdots+(2k-1)=(2k-1)k.$
Therefore,
$$k^2+k=i_1+\cdots+i_{2k+1}>(2k-1)k.$$
This inequality is not possible, if $k>2.$

Consider the case $k=2.$ We have
\begin{equation}\label{16july}
Q_{short}(V^{0,1,2},V^{0,1,2})=\lambda V^{0,1,2,3}\smile \der^0,
\end{equation}
for some $\lambda\in K.$
In obtaining this formula we use only associativity, commutativity and linearity properties of $U$
and Leibniz rule for derivation. Therefore, (\ref{16july}) is true for
any associative commutative algebra $U$ with derivation $\der$ and
$\lambda$ does not depend from $U$ and $\der.$ In particular, we can take $U=K[x]$ and $\der=\der/\der x.$ We have
  $$Q_{short}(V^{0,1,2},V^{0,1,2})(1,x,x^2,x^3,1)=\lambda V^{0,1,2,3}\smile
\der^0(1,x,x^2,x^3,1).$$
Further,
  $$Q_{short}(V^{0,1,2},V^{0,1,2})(1,x,x^2,x^3,1)=$$

$$V^{0,1,2}(1,x,V^{0,1,2}(x^2,x^3,1))
-V^{0,1,2}(1,x^2,V^{0,1,2}(x,x^3,1))
$$
$$+V^{0,1,2}(1,x^3,V^{0,1,2}(x,x^2,1))
+V^{0,1,2}(x,x^2,V^{0,1,2}(1,x^3,1))$$
$$
-V^{0,1,2}(x,x^3,V^{0,1,2}(x,x^2,1))
+V^{0,1,2}(x^2,x^3,V^{0,1,2}(1,x,1))
$$

$$=\der^2(V^{0,1,2}(1,x^2,x^3))
-V^{0,1,2}(1,x^2,6x))
+V^{0,1,2}(1,x^3,2x^0)$$
$$
+V^{0,1,2}(x,x^2,0)
-V^{0,1,2}(x,x^3,2x^0)
+V^{0,1,2}(x^2,x^3,0)
=12$$
and
$$\lambda V^{0,1,2,3}\smile \der^0(1,x,x^2,x^3,1)=\lambda 12$$
Thus, $\lambda=1.$

\begin{lm} \label{arman}
Let $\psi\in C^{k}(A,A)$ and $X$ be a linear span of the set
$\{a_1,\ldots,a_{k-1}\}$ and $Y$ be a linear span of the set
$\{b_1,\ldots,b_k\}.$ Then
$Q\,\psi(a_1,\ldots,a_{k-1},b_1,\ldots,b_k)=0,$ if
$X\subseteq Y.$
\end{lm}

{\bf Proof.} Let $X\subseteq Y$ and $dim\,Y=l\le k.$

Since $\psi$ is skew-symmetric, if $l<k,$ then
$$\psi(b_1,\ldots,b_k)=0,$$
for any $b_1,\ldots,b_k\in Y$ and
$$\psi(a_1,\ldots,a_{k-1},b_i)=0,$$
for any $a_1,\ldots,a_{k-1}\in X\subseteq Y, b_i\in Y.$ Therefore,
in this case
$Q\,\psi=0(a_1,\ldots,a_{k-1},b_1,\ldots,b_k)=0.$

Suppose now that $dim\,Y=k.$
If $dim\,X<k-1,$ by the same reasons,
$$\psi(e_{i_1},\ldots,e_{i_{k-1}},c)=0,$$
for any $c\in A.$ Thus in this case our lemma is true.

It remains to consider the case $dim\,X=k-1$ and $dim\,Y=k.$

Let $\{e_1,\ldots,e_k\}$ be basis of $Y,$ such that
$\{e_1,\ldots,e_{k-1}\}$ be basis of $X.$
Since $Q\,\psi$ is polylinear by arguments $a_1,\ldots,a_{k-1},
b_1,\ldots,b_k,$ to prove our lemma it is enough to establish,
that
$$Q\,\psi(e_{1},\ldots,e_{k-1},e_1,\ldots,e_k)=0.$$
Since $\psi$ is skew-symmetric,
$$\psi(e_{1},\ldots,e_{k-1},e_i)= 0,$$
if $i\le k-1.$
So,
$$\sum_{i=1}^k\psi(e_1,\ldots,e_{i-1},\psi(e_{1},\ldots,e_{k-1},e_i),
e_{i+1},\ldots,e_k)=$$
$$
\psi(e_1,\ldots,e_{k-1},\psi(e_{1},\ldots,e_{k-1},e_k)).$$
In other words, $Q\,\psi(e_1,\ldots,e_{k-1},e_1,\ldots,e_k)=0.$
Lemma is proved completely.

\begin{lm}\label{0123} $V^{0,1,2,3}$ is $4$-Lie multiplication,
if $p=3.$
\end{lm}

{\bf Proof.} We have
$$| V^{0,1,2,3}|=6\Rightarrow | Q\,V^{0,1,2,3}|=12.$$
The set of all $(3,4)$-partitions of $12$ is
$$\Gamma_{3,4}(12)=\left\{
(\{0,1,2\},\{0,1,2,6\}), (\{0,1,2\},\{0,1,3,5\}),\right.$$
$$ (\{0,1,2\},\{0,2,3,4\}), (\{0,1,3\},\{0,1,2,5\}),$$
$$(\{0,1,3\},\{0,1,3,4\}),
(\{0,1,4\},\{0,1,2,4\}),$$
$$(\{0,1,5\},\{0,1,2,3\}),
(\{0,2,3\},\{0,1,2,4\}), (\{0,2,4\},\{0,1,2,3\}),$$
$$\left. (\{1,2,3\},\{0,1,2,3\})
\right\}
$$
Thus,
\begin{equation}\label{arman3}
Q\,V^{0,1,2,3}=\sum_{(\alpha,\beta)\in\Gamma_{3,4}(12)}
\lambda_{(\alpha,\beta)}V^{\alpha}\smile V^{\beta},
\end{equation}
where
$$\alpha=\{i_1,i_2,i_3\}, \beta=\{i_4,i_5,i_6,i_7\},$$
$$0\le i_1<i_2<i_3, 0\le i_4<i_5<i_6<i_7, i_1+\cdots+i_7=12.$$
In receiving of formula (\ref{arman3}) one uses only the Leibniz rule.
Therefore, this formula is universal, i.e.,  coefficients
$\lambda_{\alpha,\beta}$ are integers that
do not depend on $U$ and derivation $\der.$ So, it is true for any
associative commutative algebra $U$ with derivation $\der$
In particular, we can take $U={\bf Q}[x], \der=\der/\der_x.$
To find $\lambda_{\alpha,\beta}$ we can substitute
$a_l=x^{i_l}, l=1,\ldots,7$ and calculate $Q\,V^{0,1,2,3}$ in $k[x].$
We have
$$\lambda_{\alpha,\beta}=\frac{1}{i_1!\cdots i_7!}Q\,V^{0,1,2,3}.$$

By lemma~\ref{arman},
$$Q\,V^{0,1,2,3}(1,x,x^2,1,x,x^2,x^6)=0,$$
$$Q\,V^{0,1,2,3}(1,x,x^3,1,x,x^3,x^4)=0,$$
$$Q\,V^{0,1,2,3}(1,x,x^4,1,x,x^2,x^4)=0.$$
Thus, if
$$(\alpha,\beta)\in \{(\{0,1,2\},\{0,1,2,6\}), (\{0,1,3\},\{0,1,3,4\}),
(\{0,1,4\},\{0,1,2,4\})\},$$
then $\lambda_{\alpha,\beta}=0.$

We have
$$V^{0,1,2,3}(1,x^2,x^3,x^4)=\left|\begin{array}{cccc}
1&x^2&x^3&x^4\\
0&2x&3x^2&4x^3\\
0&2&6x&12x^2\\
0&0&6&24x\\
\end{array}
\right|
=48 x^3,$$
$$V^{0,1,2,3}(1,x,x^2,x^4)=\left|\begin{array}{cccc}
1&x&x^2&x^4\\
0&1&2x&4x^3\\
0&0&2&12x^2\\
0&0&0&24x\\
\end{array}
\right|
=48 x,$$
Thus
$$V^{0,1,2,3}(1,x,x^2,1,x^2,x^3,x^4)=$$ $$
V^{0,1,2,3}(1,x,x^2,V^{0,1,2,3}(1,x^2,x^3,x^4))-0-0-
V^{0,1,2,3}(1,x^2,x^3,V^{0,1,2,3}(1,x,x^2,x^4))=$$
$$48 V^{0,1,2,3}(1,x,x^2,x^3)-48 V^{0,1,2,3}(1,x^2,x^3, x)=0,$$
and
$$\lambda_{\{0,1,2\},\{0,2,3,4\}}=0.$$

Further,
$$V^{0,1,2,3}(1,x,x^3,x^5)=\left|\begin{array}{cccc}
1&x&x^3&x^5\\
0&1&3x^2&5x^4\\
0&0&6x&20x^3\\
0&0&6&60x^2\\
\end{array}
\right|
=240x^3,$$
$$V^{0,1,2,3}(1,x,x^2,x^5)=\left|\begin{array}{cccc}
1&x&x^2&x^5\\
0&1&2x&5x^4\\
0&0&2&20x^3\\
0&0&0&60x^2\\
\end{array}
\right|
=120x^2.$$
Thus
$$V^{0,1,2,3}(1,x,x^2,1,x,x^3,x^5)=$$
$$
V^{0,1,2,3}(1,x,x^2,V^{0,1,2,3}(1,x,x^3,x^5))-0-0-
V^{0,1,2,3}(1,x,x^3,V^{0,1,2,3}(1,x,x^2,x^5))=$$
$$240 V^{0,1,2,3}(1,x,x^2,x^3)-120 V^{0,1,2,3}(1,x,x^3, x^2)=$$
$$360 V^{0,1,2,3}(1,x,x^2,x^3)=4320,$$
and
$$\lambda_{\{0,1,2\},\{0,1,3,5\}}=\frac{4320}{0!1!2!0!1!3!5!}=3.$$

Similar calculations show that
$$\lambda_{\{0,1,3\},\{0,1,2,5\}}=-3,
\lambda_{\{0,1,5\},\{0,1,2,3\}}=3,
\lambda_{\{0,2,3\},\{0,1,2,4\}}=\lambda_{\{0,2,4\},\{0,1,2,3\}}=0
$$

So, we have established that
$$
Q\,V^{0,1,2,3}=3\,V^{0,1,2}\,
V^{0,1,3,5} -3\,V^{0,1,3}\, V^{0,1,2,5} +3\,V^{0,1,5}\,V^{0,1,2,3}
$$
In particular, if $p=3,$ then
$$Q\,V^{0,1,2,3}=0.$$

\begin{crl}\label{crl0123} If $p=3,$ then $V^{1,2,3}$ is $3$-Lie multiplication and
$V^{2,3}$ is $2$-Lie multiplication.
\end{crl}

{\bf Proof.}
Follows from lemma \ref{0123} and proposition \ref{filipov1} and the
following observations
$$i(1)V^{0,1,2,3}=V^{1,2,3},$$
$$i(x)V^{1,2,3}=V^{2,3}.$$

\begin{lm}\label{p=2}
 $V^{0,1,2,3,4}$ is $5$-Lie multiplication, if $p=2.$
\end{lm}

{\bf Proof} is similar to the proof of lemma~\ref{0123}, so we will omit
calculation details. We have
$$|V^{0,1,2,3,4}|=10\Rightarrow |Q\,V^{0,1,2,3}|=20.$$
Therefore, $Q\,V^{0,1,2,3}$ is a linear combination of $V^{\alpha}\smile V^{\beta},$ where $(\alpha,\beta)\in \Gamma_{4,5}(20),$ and
$$\Gamma_{4,5}(20)=\left\{
(\{0,1,2,3\},\{0,1,2,3,8\}), (\{0,1,2,3\},\{0,1,2,4,7\}),\right.$$
$$(\{0,1,2,3\},\{0,1,2,5,6\}), (\{0,1,2,3\},\{0,1,3,4,6\}),
(\{0,1,2,3\},\{0,2,3,4,5\}),$$
$$
(\{0,1,2,4\},\{0,1,2,3,7\}),
(\{0,1,2,4\},\{0,1,2,4,6\}),
(\{0,1,2,4\},\{0,1,3,4,5\}),$$
$$
(\{0,1,2,5\},\{0,1,2,3,6\}),
(\{0,1,2,5\},\{0,1,2,4,5\}),
(\{0,1,3,4\},\{0,1,2,3,6\}),
$$
$$
(\{0,1,3,4\},\{0,1,2,4,5\}),
(\{0,1,2,6\},\{0,1,2,3,5\}),
(\{0,1,3,5\},\{0,1,2,3,5\}),
$$
$$
(\{0,2,3,4\},\{0,1,2,3,5\}),
(\{0,1,2,7\},\{0,1,2,3,4\}),
(\{0,1,3,6\},\{0,1,2,3,4\}), $$
$$
\left.(\{0,1,4,5\},\{0,1,2,3,4\}),
(\{0,2,3,5\},\{0,1,2,3,4\}),
(\{1,2,3,4\},\{0,1,2,3,4\})
\right\}
$$
Thus, there exist $\lambda_{\alpha,\beta}\in {\bf Z},$ such that
$$
Q\,V^{0,1,2,3,4}=\sum_{(\alpha,\beta)\in\Gamma_{4,5}(20)}
\lambda_{(\alpha,\beta)}V^{\alpha}\smile V^{\beta},
$$

Calculations as in the proof of lemma \ref{0123} show that
$$Q\,V^{0,1,2,3,4} =$$
$$
  4\,V^{0, 1, 2, 7}\,V^{0, 1, 2, 3, 4} +
  2\,V^{0, 1, 3, 6}\,V^{0, 1, 2, 3, 4} -
  2\,V^{0, 1, 4, 5}\,V^{0, 1, 2, 3, 4}$$
  $$ +
  2\,V^{0, 2, 3, 5}\,V^{0, 1, 2, 3, 4} +
  2\,V^{0, 1, 2, 6}\,V^{0, 1, 2, 3, 5} -
  2\,V^{0, 2, 3, 4}\,V^{0, 1, 2, 3, 5}$$
  $$ -
  2\,V^{0, 1, 2, 5}\,V^{0, 1, 2, 3, 6} -
  2\,V^{0, 1, 3, 4}\,V^{0, 1, 2, 3, 6} -
  4\,V^{0, 1, 2, 4}\,V^{0, 1, 2, 3, 7} $$
  $$+
  2\,V^{0, 1, 3, 4}\,V^{0, 1, 2, 4, 5} +
  4\,V^{0, 1, 2, 3}\,V^{0, 1, 2, 4, 7} +
  2\,V^{0, 1, 2, 3}\,V^{0, 1, 2, 5, 6}$$
  $$ -
  2\,V^{0, 1, 2, 4}\,V^{0, 1, 3, 4, 5} +
  2\,V^{0, 1, 2, 3}\,V^{0, 1, 3, 4, 6} +
  2\,V^{0, 1, 2, 3}\,V^{0, 2, 3, 4, 5}
$$
In particular, if $p=2,$ then
$$Q\,V^{0,1,2,3,4}=0.
$$

\begin{crl} \label{crl01234} If  $p=2,$ then
$V^{1,2,3,4}$ is $4$-Lie multiplication,
$V^{2,3,4}$ is $3$-Lie multiplication and
$V^{3,4}$ is $2$-Lie multiplication.
\end{crl}

{\bf Proof.} Follows from lemma \ref{p=2} and proposition
\ref{filipov1} and the following facts
$$V^{1,2,3,4}=i(1)V^{0,1,2,3,4},$$
$$V^{2,3,4}=i(x)V^{1,2,3,4},$$
$$V^{3,4}=i(x^2)V^{2,3,4}/2.$$

{\bf Remark.} It is easy to calculate $Q\,V^{1,2,3,4},
Q\,V^{2,3,4}, Q\,V^{3,4}$ and $Q\,V^{1,2,3}, Q\,V^{2,3}$
over ${\bf Z}$ directly as in calucating $Q\,V^{0,1,2,3,4}$
and $Q\,V^{0,1,2,3}.$ For example,
$$Q\,V^{1,2,3,4}=$$ $$
  4\,V^{1, 2, 7}\,V^{1, 2, 3, 4} +
  2\,V^{1, 3, 6}\,V^{1, 2, 3, 4} -
  2\,V^{1, 4, 5}\,V^{1, 2, 3, 4} $$
$$+
  2\,V^{2, 3, 5}\,V^{1, 2, 3, 4} +
  2\,V^{1, 2, 6}\,V^{1, 2, 3, 5} -
  2\,V^{2, 3, 4}\,V^{1, 2, 3, 5} $$
  $$-
  2\,V^{1, 2, 5}\,V^{1, 2, 3, 6} -
  2\,V^{1, 3, 4}\,V^{1, 2, 3, 6} -
  4\,V^{1, 2, 4}\,V^{1, 2, 3, 7} $$
  $$ +
  2\,V^{1, 3, 4}\,V^{1, 2, 4, 5} +
  4\,V^{1, 2, 3}\,V^{1, 2, 4, 7} +
  2\,V^{1, 2, 3}\,V^{1, 2, 5, 6} $$
  $$-
  2\,V^{1, 2, 4}\,V^{1, 3, 4, 5} +
  2\,V^{1, 2, 3}\,V^{1, 3, 4, 6} +
  2\,V^{1, 2, 3}\,V^{2, 3, 4, 5}
$$
and
$$Q\,V^{1,2,3}=$$
$$3\,V^{1,2}\,V^{1,3,5}-3\,V^{1,3}\,V^{1,2,5}+3\,V^{1,5}\,V^{1,2,3}.$$

\begin{lm} ({\rm Lukas}) If $a=\sum_{i\ge 0}a_ip^i,$ $ 0\le a_i<p,$ and
$b=\sum_{j\ge 0}b_jp^j,$ $0\le b_i<p,$ are $p$-adic presentations of
$a$ and $b,$ then
$${a\choose b}\equiv \prod_i{a_i\choose b_i}(mod\,p)$$
\end{lm}

Let
$$O_1(m)= \{x^{(i)}: x^{(i)}x^{(j)}={i+j\choose i}x^{(i+j)}, 0\le i,j<p^m\}$$
be divided power algebra over the field $K$ of characteristic $p>0$ and
$$\der: O_1(m)\rightarrow O_1(m), \quad x^{(i)}\mapsto x^{(i-1)},$$
be special derivation.

\begin{lm}\label{one} Let $U=O_1(m), p>0.$ Then
\begin{enumerate}
\renewcommand{\theenumi}{\roman{enumi}}
\item $\der^q\in Der\,U$ if and only if $q=p^k$ for some $k>0.$
\item $\der^{p^k-1}\wedge \der^{p^k}$ is $2$-Lie, if and only if
$p=2$ or $p=3, k=1.$
\item $\der^{p^k-2}\wedge \der^{p^k-1}\wedge\der^{p^k}$ is
$3$-Lie, if and only if $p=3, k=1,$ or $p=2, k=1$ or $p=2,k=2.$
\end{enumerate}
\end{lm}

{\bf Proof.} {\rm (i)} If $\der^q$ derivation, then
$$\der^q(x^{(i)}x^{(j)})=\der^q(x^{(i)})x^{(j)}+x^{(i)}\der^q(x^{(j)})\Rightarrow$$
$${i+j\choose i}-{i+j-q\choose i}-{i+j-q\choose j}\equiv
0(mod\,p),$$ for any $0\le i,j.$ If $q$ is not power of $p,$ say
$q=p^kb,$ where $b>1$ and $p$ are mutually prime, then by Lukas
lemma, $$i=p^k, j=p^k(b-1)\Rightarrow {p^kb\choose p^k}\equiv
b\not\equiv 0(mod\,p).$$ Contradiction. So, $\der^q\not\in
Der\,U,$ if $q$ is not power of $p.$

\medskip

{\rm (ii)} Suppose that $\psi=\der^{p^k-1}\wedge\der^{p^k}$ is
$2$-Lie. Then
$$A(i,j,s):=\psi(x^{(i)},\psi(x^{(j)},x^{(s)}))-\psi(\psi(x^{(i)},x^{(j)}),x^{(s)})
-\psi(x^{(j)},\psi(x^{(i)},x^{(s)}))=0,$$ for any $i,j,s\ge 0.$
Take $i=p^k, j= p^k+1, s=2 p^k-3.$ We have
$$\psi(x^{(i)},x^{(j)})=\lambda x^{(2)}, \psi(x^{(i)},x^{(s)})=\mu
x^{(p^k-2)},$$ for some $\lambda,\mu\in K$ and
$$\psi(x^{(i)},\psi(x^{(j)},x^{(s)}))= \psi(x^{(p^k)},x^{(2)}
x^{(p^k-3)}-x x^{(p^k-2)})=$$ $$2\psi(x^{(p^k)},x^{(p^k-1)})=-2.$$
If $p^k>3,$ then $2<p^k-1.$ Since $\psi(x^{(i)}, u)=0,$ for any
$u\in U, i<p^k-1,$ we have $$A(p^k,p^k+1,2 p^k-3)=-2\ne 0,$$ if
$p^k>3, p\ne 2.$

Now check that $\psi$ is $2$-Lie, if $(p,k)=(3,1)$ or $p=2.$

\medskip

Let $p^k=3,$ i.e., $p=3, k=1.$ Prove that $\der^2\wedge\der^3$ is
$2$-Lie operation on $U.$ Let
$$\omega(a,b)=\der^2(a)\der^3(b)-\der^3(a)\der^2(b).$$
We have
$$\omega(\omega(a,b),c)=$$

$$ \der^2(\der^2(a)\der^3(b))\der^3(c)
-\der^2(\der^3(a)\der^2(b))\der^3(c)
-\der^3(\der^2(a)\der^3(b))\der^2(c)
+\der^3(\der^3(a)\der^2(b))\der^2(c)=$$

$$\mathop{\der^4(a)\der^3(b))\der^3(c)}\limsim+
\mathop{2\der^3(a)\der^4(b)\der^3(c)}\limline+
\der^2(a)\der^5(b)\der^3(c)$$ $$
 -\der^5(a)\der^2(b)\der^3(c)-\mathop{2\der^4(a)\der^3(b)\der^3(c)}\limsim
 -\mathop{\der^3(a)\der^4(b)\der^3(c)}\limline$$
$$ -\der^5(a)\der^3(b)\der^2(c)-\der^2(a)\der^6(b)\der^2(c)
+\der^6(a)\der^2(b)\der^2(c)+\der^3(a)\der^5(b)\der^2(c)=$$

$$-\der^4(a)\der^3(b))\der^3(c)+\der^3(a)\der^4(b)\der^3(c)+
\der^2(a)\der^5(b)\der^3(c) -\der^5(a)\der^2(b)\der^3(c)$$ $$
-\der^5(a)\der^3(b)\der^2(c)-\der^2(a)\der^6(b)\der^2(c)
+\der^6(a)\der^2(b)\der^2(c)+\der^3(a)\der^5(b)\der^2(c)$$

By similar way one can easy find  $\omega(\omega(b,c),a)$ and
$\omega(\omega(c,a),b)$. Therefore,

$$\omega(\omega(a,b),c)+ \omega(\omega(b,c),a)+
\omega(\omega(c,a),b)=$$

$$-\mathop{\der^4(a)\der^3(b))\der^3(c)}\limline+
\mathop{\der^3(a)\der^4(b)\der^3(c)}\limsim +
\mathop{\der^2(a)\der^5(b)\der^3(c)}\limsimeq
-\mathop{\der^5(a)\der^2(b)\der^3(c)}\limeq $$

$$ -\mathop{\der^5(a)\der^3(b)\der^2(c)}\limequiv-
\mathop{\der^2(a)\der^6(b)\der^2(c)}\limsmile
+\mathop{\der^6(a)\der^2(b)\der^2(c)}\limapprox\
+\mathop{\der^3(a)\der^5(b)\der^2(c)}\limcong$$

$$-\mathop{\der^4(b)\der^3(c))\der^3(a)}\limsim
+\mathop{\der^3(b)\der^4(c)\der^3(a)}\limasymp+
\mathop{\der^2(b)\der^5(c)\der^3(a)}\limdoteq
-\mathop{\der^5(b)\der^2(c)\der^3(a)}\limcong$$

$$ -\mathop{\der^5(b)\der^3(c)\der^2(a)}\limsimeq
-\mathop{\der^2(b)\der^6(c)\der^2(a)}\limfrown
+\mathop{\der^6(b)\der^2(c)\der^2(a)}\limsmile
+\mathop{\der^3(b)\der^5(c)\der^2(a)}\limli$$

$$-\mathop{\der^4(c)\der^3(a))\der^3(b)}\limasymp
+\mathop{\der^3(c)\der^4(a)\der^3(b)}\limline+
\mathop{\der^2(c)\der^5(a)\der^3(b)}\limequiv-
 \mathop{\der^5(c)\der^2(a)\der^3(b)}\limli$$ $$
-\mathop{\der^5(c)\der^3(a)\der^2(b)}\limdoteq
-\mathop{\der^2(c)\der^6(a)\der^2(b)}\limapprox
+\mathop{\der^6(c)\der^2(a)\der^2(b)}\limfrown
+\mathop{\der^3(c)\der^5(a)\der^2(b)}\limeq$$ $$=0.$$

\medskip

Let now $p=2.$ Prove that for any nonnegative integer $k,$
$\der^{2^k}\wedge \der^{2^k-1}$ is $2$-Lie operation on $U.$ For
$$\omega=\der^{2^k}\wedge \der^{2^k-1}$$ we have
$$\omega(a,b)=\der(\der^{2^k-1}(a)\der^{2^k-1}(b)).$$

Therefore $$\omega(\omega(a,b),c)=$$
$$\der(\der^{2^k}(\der^{2^k-1}(a)\der^{2^k-1}(b))\der^{2^k-1}(c))=$$
$$\der\left(\der^{2^{k+1}-1}(a)\der^{2^k-1}(b)\der^{2^k-1}(c)
+\der^{2^k-1}(a)\der^{2^{k+1}-1}(b)\der^{2^k-1}(c)\right) =$$

Similarly, $$\omega(\omega(b,c),a)=$$
$$\der\left(\der^{2^{k+1}-1}(b)\der^{2^k-1}(c)\der^{2^k-1}(a)
+\der^{2^k-1}(b)\der^{2^{k+1}-1}(c)\der^{2^k-1}(a)\right),$$ and
$$\omega(\omega(c,a),b)=$$
$$\der\left(\der^{2^{k+1}-1}(c)\der^{2^k-1}(a)\der^{2^k-1}(b)
+\der^{2^k-1}(c)\der^{2^{k+1}-1}(a)\der^{2^k-1}(b)\right).$$ Thus,
$$\omega(\omega(a,b),c)+ \omega(\omega(b,c),a)+
\omega(\omega(c,a),b)=$$

$$\der\left(\mathop{\der^{2^{k+1}-1}(a)\der^{2^k-1}(b)\der^{2^k-1}(c)}\limline
+\mathop{\der^{2^k-1}(a)\der^{2^{k+1}-1}(b)\der^{2^k-1}(c)}\limsim\right.$$
$$+\mathop{\der^{2^{k+1}-1}(b)\der^{2^k-1}(c)\der^{2^k-1}(a)}\limsim
+\mathop{\der^{2^k-1}(b)\der^{2^{k+1}-1}(c)\der^{2^k-1}(a)}\limeq$$
$$\left.+\mathop{\der^{2^{k+1}-1}(c)\der^{2^k-1}(a)\der^{2^k-1}(b)}\limeq
+\mathop{\der^{2^k-1}(c)\der^{2^{k+1}-1}(a)\der^{2^k-1}(b)}\limline\right)$$
$$=0.$$

\medskip

{\rm (iii)} Suppose that $\psi=\der^{2^k-2}\wedge\der^{2^k-1}\wedge
\der^{2^k}, k>2,$ is $3$-Lie.

Then $$A(i_1,i_2,i_3,i_4,i_5)=
B(i_1,\ldots,i_5)-C(i_1,\ldots,i_5)=0,$$ for any
$i_1,\ldots,i_5\ge 0,$ where $$B(i_1,\ldots,i_5)=
\psi(x^{(i_1)},x^{(i_2)},\psi(x^{(i_3)},x^{(i_4)},x^{(i_5)})), $$
$$C(i_1,\ldots,i_5)=$$
$$
\psi(\psi(x^{(i_1)},x^{(i_2)},x^{(i_3)}),x^{(i_4)},x^{(i_5)})+
\psi(x^{(i_3)},\psi(x^{(i_1)},x^{(i_2)},x^{(i_4)}),x^{(i_5)}+$$
$$\psi(x^{(i_3)},x^{(i_4)},\psi(x^{(i_1)},x^{(i_2)},x^{(i_5)})).$$
Take $i_1=2^k-2, i_2=2^k, i_3=2^k-1, i_4=2^k+1, i_5=2^{k+1}-4.$
Then $$\psi(x^{(i_3)},x^{(i_4)},x^{(i_5)})= \left|
\begin{array}{ccc}
x^{(1)}&x^{(3)}&x^{(2^{k}-2)}\\ x^{(0)}&x^{(2)}&x^{(2^{k}-3)}\\
0&x^{(1)}&x^{(2^{k}-4)}\\
\end{array}\right|= x^{(2^k-1)}$$
and $$\psi(x^{(i_1)},x^{(i_2)},x^{(i_3)})=\lambda x^{(0)},$$
 $$\psi(x^{(i_1)},x^{(i_2)},x^{(i_4)})=\mu x^{(2)},$$
 $$\psi(x^{(i_1)},x^{(i_2)},x^{(i_5)})=\nu x^{(2^k-3)},$$
for some $\lambda,\mu,\nu\in K.$ If $k>2,$ then $2<2^k-2.$
Therefore, $$B(i_1,\ldots,i_5)=1,$$ $$C(i_1,\ldots,i_5)=0,$$ if
$k>2.$ So, if $k>2$ we obtain contradiction $$A(2^k-2, 2^k, 2^k-1,
2^k+1,2^{k+1}-4)=1\ne 0.$$

It is easy to check that $\der^2\wedge\der^3\wedge\der^4$ and
$\der^0\wedge\der\wedge \der^2$ are  $3$-Lie.

\begin{crl}{\rm (i)}$p=3.$ For any  $k\in {\bf Z}_+,$
$\der^{p^k-p^{k-1}}\wedge \der^{p^k}$ is 2-Lie operation.

{\rm (ii)} $p=2.$ For any $k,l\in {\bf Z}_+,  k>l,$
$\der^{p^k-p^l}\wedge\der^{p^k}$ is 2-Lie operation.
\end{crl}

{\bf Proof.} i) For $p=3,$ we have $p^k-p^{k-1}=2 p^{k-1},
p^k=3p^{k-1}.$ Since $F=\der^{p^{k-1}}\in Der\,U,$ our statement
follows from lemma \ref{one} {\rm ii} used for $F$ instead of $\der.$

ii) For $p=2,$ we have $p^k-p^l=p^l(p^{k-l}-1), p^k=p^{k-l}p^l.$
Therefore, for $F=\der^{p^l},$
$$\der^{p^k}=F^{p^{k-l}},\der^{p^k-p^l}=F^{p^{k-l}-1}.$$ Our
statement follows from lemma \ref{one} {\rm ii} used for $F^{p^{k-l}}$ instead
of $\der^{p^k}.$

{\bf Proof of theorem \ref{trieste1}}

{\rm i} Corollary \ref{altynaitugankun}

{\rm ii} Lemma \ref{altynkustugankun}

{\rm iii} Suppose that $(U,V^{0,1,\ldots,q})$ is $(q+1)$-Lie.
If $q=1,$ then it is $2$-Lie for any characteristic $p.$.

Assume that $q>1.$  By lemma \ref{filipov1},
$V^{q}=i(1)i(x)\cdots i(x^{(q-2)}) V^{0,1,\ldots,q}$ is $1$-Lie, i.e.,
$\der^q\in Der\,U.$ If $q>1$ and $p=0$ this is not possible.

So $p>0.$ Take $U=O_1(m).$ By lemma \ref{one} {\rm i},
$q$ should be power of $p.$ Suppose that $k=p^t.$

By lemma \ref{filipov1}, $V^{p^t-1,p^t}=i(1)i(x)\cdots i(x^{(p^t-2)})
V^{0,1,\ldots,p^t}$ is $2$-Lie. By lemma \ref{one} {\rm ii}
it is possible in the following cases $p=2$ or $p=3,t=1.$

By lemma \ref{filipov1}, $V^{2^t-2,2^t-1,2^t}=
i(1)i(x)\cdots i(x^{(2^t-3)})V^{0,1,\ldots,2^t}$ is $3$-Lie.
By lemma \ref{one} {\rm iii} this is possible only in the cases
$p=2,t=1$ or $p=2, t=2.$
Theorem \ref{trieste1} is proved completely.

{\bf Proof of theorem \ref{trieste2}}
Corollaries \ref{crl0123} and \ref{crl01234}. Calculations for
checking that
$\sum_{i=1}^{2^l}V^{0,i,2^l+1-i}, p=2,$ is $3$-Lie can be done by analogous way.
Then by proposition \ref{filipov1} $\sum_{i=1}^{2^l}V^{i,2^l+1-i}, p=2,$ is
$2$-Lie.


\begin{thebibliography}{10}

\bibitem{Awata} H. Awata, M. Li, D. Minic, T. Yoneya,
{\em On the quantization of Nambu brackets,} hep-th/9906248 .

\bibitem{Dzhumavestnik} A.S. Dzhumadil'daev,
{\em A remark on the space of invariant differential operators},
Vestnik Moskov. Univ, Ser.1 Mat.Mekh, 1982, No.2, 49-54, 116=engl.
transl.  Moscow Univ. Math. Bull., {\bf 37}(1982), No.2, p.63-68.

\bibitem{dzhumafunctanalis}
A.S. Dzhumadil'daev,
{\em Integral and $mod p-$cohomologies of the Lie
algebra $W_1$}, Funct Anal. Pril. {\bf 22}(1988), no.3, 68-70=engl.transl.
Funct.Anal.Appl., {\bf 22}(1988), No.3, p.226-228 (1989).

\bibitem{Filipov1}
V.T. Filipov, {\em $n$-Lie algebras,} Sibirsk. Mat.Zh. {\bf 26}(1985), no. 6, 126-140, 191=engl. transl. Siberian Math. J., {\bf 26}(1985), no.6, 879-891


\bibitem{Filipov2} V.T. Filipov, {\em On $n$-Lie algebra of
jacobians,} Sibirsk.Mat.Zh., {\bf 39}(1998), No.3, 660-669=
engl. transl. Siberian Math. J., {\bf 39}(1998), no.3, 573-581.

\bibitem{Hanlon} P. Hanlon, M. Wachs, {\em On Lie $k$-algebras,}
Adv. Math., {\bf 113}(1995), 206-236.

\bibitem{Hoppe} J. Hoppe, {\em On $M$-algebras, the quantization
of Nambu-Mechanics, and the volume preserving diffeomorphisms,}
hep-th/9602020 .

\bibitem{Kurosh} A.G. Kurosh, {\em Multi-operator rings and
algebras,} Uspechi Matem.Nauk,  {\bf 24}(1969), No.1, 3-15.

\bibitem{Nambu} Y.Nambu, {\em Generalized Hamiltonian mechanics,}
Phys. Rev., D 7, 2405-2412, 1973.

\bibitem{Takh} L.A. Takhtajan, {\em On the foundation of the
generalized Nambu mechanics,} Commun. Math. Phys., {\bf
160}(1994), 295-315.

\bibitem{StashefLada}J. Stasheff, T. Lada,
{\em Introduction to SH Lie algebras fo physicists,}
Inter. J. Theor.Physics, {\bf 32}(1993), No.7, 1087-1103.
\end{thebibliography}
\end{document}